\documentclass[11pt,a4paper]{glasnik}
\usepackage{glasnik}
\usepackage{graphicx, url}
\usepackage[none]{hyphenat}
\def\CurrentVolume{}
\def\CurrentYear{}

\newarticle

\begin{document}

\title[On the multiplicity in Pillai's problem]{On the multiplicity in Pillai's problem with Fibonacci numbers and powers of a fixed prime}
\UniCountry{Makerere University, Uganda and Wits University, South Africa}
\author[H Batte$^{1,*} $, M Ddamulira$^{1}$, J Kasozi$^{1}$, and F Luca$^{2}$]{Herbert Batte$^{1,*} $, Mahadi Ddamulira$^{1}$, Juma Kasozi$^{1}$, and Florian Luca$^{2}$}
\address[F. author]{Department of Mathematics\\ Makerere University\\ Kampala\\ Uganda}
\email{hbatte91@gmail.com}
\address[S. author]{Department of Mathematics\\ Makerere University\\ Kampala\\ Uganda}
\email{mahadi.ddamulira@mak.ac.ug}
\address[T. author]{Department of Mathematics\\ Makerere University\\ Kampala\\ Uganda}
\email{juma.kasozi@mak.ac.ug}
\address[L. author]{School of Mathematics, University of the Witwatersrand, Johannesburg, South Africa; Research Group in Algebraic Structures and Applications, King Abdulaziz University, Jeddah, Saudi Arabia, Max Planck Institute for Software Systems, Saarbr\"ucken, Germany and Centro de Ciencias Matem\'aticas UNAM, Morelia, Mexico}
\email{Florian.Luca@wits.ac.za}


\noindent 
\thanks{$ ^{*} $ Corresponding author}

\abstract{ Let $ \{F_n\}_{n\geq 0} $ be the sequence of Fibonacci numbers and let $p$ be a prime. For an integer $c$ we write $m_{F,p}(c)$ for the number of distinct representations of $c$ as $F_k-p^\ell$ with $k\ge 2$ and $\ell\ge 0$. We prove that $m_{F,p}(c)\le 4$. }

\maketitle

\section{Introduction}
\label{intro}
\subsection{Background}
\label{sec:1.1}
\noindent Let $\{F_n\}_{n\ge 0}$ be the Fibonacci sequence given by $F_0=0$, $F_1=1$ and $F_{n+2}=F_{n+1}+F_{n}$ for all $n \geq 0$.
The first few terms of this sequence are given by 
$$
0,~1,~1,~2,~3,~5,~8,~13,~21,~34,~55,~89,~144,\ldots. 
$$
For fixed integers $a>1,b>1,c$, the Diophantine equation
\begin{equation} \label{1.1}
a^{x}-b^{y}=c,
\end{equation}
in nonnegative integers $x$, $y$ is known as the Pillai equation, see \cite{PIL1}. Pillai was interested if the above equation can have more than one solution $(x,y)$ and proved that if $a$ and $b$ are positive and coprime and $|c|>c_0(a,b)$, then the above equation has at most one solution $(x,y)$. Variants of the Pillai problem have been recently considered in which one takes $a$ to be $2$ or $3$ but replaces the sequence of powers 
of $b$ by some other sequence of positive integers of exponential growth such as Fibonacci numbers, Tribonacci numbers, Pell numbers and even $k$--generalized Fibonacci numbers where $k$ is also unknown. In all these works, it was shown that the conclusion of the original Pillai problem is retained (so, every large integer has at most one such representation) except for some cases where parametric families exist which are completely classified together with the  small exceptional cases with multiple such representations which have also been computed. See, for example,  \cite{BFY}, \cite{CPZ}, \cite{DFR} and \cite{DDA}.
Here, we retain the Fibonacci sequence but replace powers of $2$ or $3$ by powers of an arbitrary but fixed prime $p$. Write
$$
m_{F,p}(c):=\#\{(k,\ell): k\ge 2,~\ell\ge 0,~c=F_k-p^\ell\}.
$$
We imposed the condition $k\ge 2$ above because $F_1=F_2=1$. Our result is the following.

\subsection{Main Result}
\label{sec:1.2}
\begin{theorem}\label{4.1a} The inequality $m_{F,p}(c)\le 4$ holds for all primes $p$ and all integers $c$. 
\end{theorem}
We believe that the better result $m_{F,p}(c)\le 3$ holds but we did not succeed in proving this. Further, quite possibly  $m_{F,p}(c)=3$ holds only for finitely many pairs $(p,c)$, and maybe only for the 
following 3 pairs
 \begin{equation} \label{2b}
 \begin{matrix} 
 (2,-3): & \quad & -3  & =  & F_{7}-2^{4} & = & F_{5}-2^{3} & = & F_{2}-2^{2},\\          
 (2,0):  & \quad &  0 &  = &  F_{6}-2^{3} & = &  F_{3}-2^{1} & = &  F_{2}-2^{0},\\
 (2,1): & \quad  & 1 & = &  F_{5}-2^{2} & = &  F_{4}-2^{1} & = & F_{3}-2^{0}.
 \end{matrix}
 \end{equation}
We leave proving that $m_{F,p}(c)\le 3$ and classifying the pairs of integers $(p,c)$ with $m_{F,p}(c)=3$ as a problem to the reader. On the other hand, we believe that $m_{F,p}(c)=2$ holds for infinitely many pairs $(c,p)$. For that, it suffices to look for $c$ with two representations of the form $c=F_{k_1}-p=F_{k_2}-1$, so
 the two representations $(k_1,\ell_1)$ and $(k_2,\ell_2)$ have $\ell_1=1,~\ell_2=0$. Then $p=F_{k_1}-F_{k_2}+1$. A calculation revealed $2161$ primes $p$ of the above form in the range $2\le k_2< k_1\le 1000$. 

Before embarking to the proof, let us remark that it is not surprising that $m_{F,p}(c)$ is bounded by  an absolute constant. Indeed, letting $m_{F,p}(c)=m$, then the equation 
 $$
 c=F_k-p^{\ell}
 $$
 has $m$ solutions $(k,\ell)$ with $k\ge 2$ and $\ell\ge 0$. If $c=0$, we get $F_k=p^{\ell}$. The only solution of this equation with $\ell\ge 2$ is $F_6=2^3$ by the well-known result concerning perfect powers in the Fibonacci sequence \cite{BMS}. Thus, $m_{F,p}(0)\le 3$, and $m_{F,p}(0)=3$ holds only for the prime $p=2$. When $c\ne 0$, then by using the Binet formula for the Fibonacci sequence 
 \begin{equation}
 \label{eq:Binet}
 F_k=\frac{\alpha^k-\beta^k}{\sqrt{5}},\quad {\text{\rm where}}\quad (\alpha,\beta)=\left(\frac{1+{\sqrt{5}}}{2},\frac{1-{\sqrt{5}}}{2}\right)
 \end{equation}
 valid for all $k\ge 0$, our equation can be rewritten as
 $$
 1=\frac{1}{c{\sqrt{5}}}\alpha^k-\frac{1}{c{\sqrt{5}}}\beta^k-\frac{1}{c} p^{\ell}.
 $$
 This is a particular case of the equation $1=a_1x_1+\cdots+a_sx_s$ with $s=3$, $(a_1,a_2,\ldots,a_s)=(1/(c{\sqrt{5}}),~\pm 1/(c{\sqrt{5}}),-1/c)$ and $x_1,x_2,\ldots,x_s$ unknowns in the multiplicative group $\Gamma$ generated by $\{\alpha,p\}$ inside ${\overline{\mathbb C}}^*$ of rank $r=2$ (note that $\beta=-\alpha^{-1}$). Furthermore, such a solution is nondegenerate in the sense that no subsum $\sum_{i\in I} a_ix_i$ vanishes for some subset $I\subseteq \{1,2,\ldots,s\}$ since for us $k\ge 2$.  Theorem 6.1 in \cite{AV} immediately gives that for a fixed choice of nonzero coefficients $(a_1,a_2,\ldots,a_s)$ the number of such solutions is 
 $$
 \le (8s)^{4s^4(s+r+1)},
 $$ 
 so for us 
 $$
 m\le 2\cdot (8\cdot 3)^{4\cdot 3^4(3+2+1)}
 $$
 and the right--hand side above exceeds $10^{2500}$. Thus, while the boundedness of $m$ follows easily from known results on the finiteness of non-degenerate solutions to ${\mathcal S}$-unit equations, 
 the merit of our paper is to give a bound on $m$ which is quite close to the best possible.
\section{Methods}
We use three times Baker-type lower bounds for nonzero linear forms in two or three logarithms of algebraic numbers. There are many such bounds mentioned in the literature like that of Baker and W{\"u}stholz from \cite{BW} or Matveev from \cite{MAT}. Before we can formulate such inequalities we need the notion of height of an algebraic number recalled below.

\begin{definition}\label{3.1}
Let $ \gamma $ be an algebraic number of degree $ d $ with minimal primitive polynomial over the integers $$ a_{0}x^{d}+a_{1}x^{d-1}+\cdots+a_{d}=a_{0}\prod_{i=1}^{d}(x-\gamma^{(i)}), $$ where the leading coefficient $ a_{0} $ is positive. Then, the logarithmic height of $ \gamma$ is given by $$ h(\gamma):= \dfrac{1}{d}\Big(\log a_{0}+\sum_{i=1}^{d}\log \max\{|\gamma^{(i)}|,1\} \Big). $$
\end{definition}
In particular, if $ \gamma$ is a rational number represented as $\gamma=p/q$ with coprime integers $p$ and $ q\ge 1$, then $ h(\gamma ) = \log \max\{|p|, q\} $. 
\newpage
The following properties of the logarithmic height function $ h(\cdot) $ will be used in the rest of the paper without further reference:
\begin{equation}\nonumber
\begin{aligned}
h(\gamma_{1}\pm\gamma_{2}) &\leq h(\gamma_{1})+h(\gamma_{2})+\log 2;\\
h(\gamma_{1}\gamma_{2}^{\pm 1} ) &\leq h(\gamma_{1})+h(\gamma_{2});\\
h(\gamma^{s}) &= |s|h(\gamma)  \quad {\text{\rm valid for}}\quad s\in \mathbb{Z}.
\end{aligned}
\end{equation}
A linear form in logarithms is an expression
\begin{equation}
\label{eq:Lambda}
\Lambda:=b_1\log \alpha_1+\cdots+b_t\log \alpha_t,
\end{equation}
where for us $\alpha_1,\ldots,\alpha_t$ are positive real  algebraic numbers and $b_1,\ldots,b_t$ are nonzero integers. We assume, $\Lambda\ne 0$. We need lower bounds 
for $|\Lambda|$. We write ${\mathbb K}={\mathbb Q}(\alpha_1,\ldots,\alpha_t)$ and $D$ for the degree of ${\mathbb K}$.
We start with the general form due to Matveev \cite{MAT}. 

\begin{theorem}
\label{thm:Mat} 
Put $\Gamma:=\alpha_1^{b_1}\cdots \alpha_t^{b_t}-1=e^{\Lambda}-1$. Then 
$$
\log |\Gamma|>-1.4\cdot 30^{t+3}\cdot t^{4.5} \cdot D^2 (1+\log D)(1+\log B)A_1\cdots A_t,
$$
where $B\ge \max\{|b_1|,\ldots,|b_t|\}$ and $A_i\ge \max\{Dh(\alpha_i),|\log \alpha_i|,0.16\}$ for $i=1,\ldots,t$.
\end{theorem}
We continue with $t=2$. Let $A_1>1,~A_2>1$ be real numbers such that 
\begin{equation}
\label{eq:Ai}
\log A_i\ge \max\left\{h(\alpha_i),\frac{|\log \alpha_i|}{D},\frac{1}{D}\right\}\quad {\text{\rm for}}\quad i=1,2.
\end{equation}
Put
$$
b':=\frac{|b_1|}{D\log A_2}+\frac{|b_2|}{D\log A_1}.
$$
The following  result is Corollary 2 in \cite{LMN}. 
\begin{theorem}
\label{thm:LMN}
In case $t=2$, we have  
$$
\log |\Lambda|\ge -24.34 D^4\left(\max\left\{\log b'+0.14,\frac{21}{D},\frac{1}{2}\right\}\right)^2\log A_1\log A_2.
$$
\end{theorem}

After some calculations with the above theorems, we end up with some upper bounds on our variables which are too large, thus we need to reduce them. We use the following 
result of Legendre which is related to continued fractions (see Theorem 8.2.4 in \cite{ME}).

\begin{lemma} 
\label{lem:Legendre}
Let $ \tau $ be an irrational number with continued fraction $[a_0,a_1,\ldots]$ and convergents ${p_{0}}/{q_{0}}$, ${p_{1}}/{q_{1}}, \ldots $. Let $M$ be a positive integer. Let $ N $ be a non-negative integer such that
	$ q_{N} > M $.\\
	Then putting $ a(M) := \max \{a_{i}: 0\le i\le N   \}$, the inequality
	$$ \bigg| \tau-\frac{r}{s}   \bigg| > \dfrac{1}{(a(M)+2)s^2},  $$
	holds for all pairs $ (r, s) $ of positive integers with $ 0 < s < M $.  
\end{lemma}

Let $\{L_n\}_{n\ge 0}$ be the Lucas sequence given by $L_0=2,~L_1=1$ and $L_{n+2}=L_{n+1}+L_n$ for all $n\ge 0$. We need the following lemma.

\begin{lemma}
\label{lem:FL}
If $m\ge n$ and $m\equiv n\pmod 2$, then 
$$
F_m-F_n=F_{(m-\delta n)/2}L_{(m+\delta n)/2},\quad {\text{\rm where}}\quad \delta=(-1)^{(m-n)/2}.
$$
\end{lemma}

For  a prime $p$ let $z(p)$ be {\it the order of appearance} of the prime $p$ in the Fibonacci sequence (sometimes also called {\it the entry point} of $p$) which is the smallest positive integer $k$ such that 
$p\mid F_k$. This exists for every prime number $p$. It is well-known that $p\mid F_k$ if and only if $z(p)\mid k$. Further, writing $\nu_p(m)$ for the exponent of $p$ in the factorization of $m$, and putting $e_p:=\nu_p(F_{z(p)})$, then it is well-known that whenever $p\mid F_k$ we have $\nu_p(F_k)\ge e_p$ and further if $f=\nu_p(F_k)$ is positive then $p^{f-e_p}z(p)\mid k$. The following computational result is due to McIntosh and Roettger \cite{MR}.

\begin{lemma}
\label{lem:MR}
If $p<10^{14}$, then $p\| F_{z(p)}$.
\end{lemma}  

Finally, we present an analytic argument which is Lemma 7 from \cite{GL}. It is useful when obtaining upper bounds on some positive real variable involving powers of the logarithm of the variable itself.
\begin{lemma}\label{3.3b} If $ s \geq 1 $, $T > (4s^2)^s$ and $T > \displaystyle \frac{x}{(\log x)^s}$, then $$x < 2^s T (\log T)^s.$$	
\end{lemma}
In the addition to the above results, we also used computations with Mathematica.

\section{Proof of Theorem \ref{4.1a}}
\subsection{Notation}
\label{sec:3.1}

From now on, we work with pairs $(p,c)$ such that $m_{F,p}(c)\ge 3$.  We may assume that $p\ge 5$ since the cases $p=2,3$ were treated in \cite{DFR} and \cite{DDA}, respectively. 
We write $m:=m_{F,p}(c)$ and write 
\begin{equation}
\label{eq:41}
c=F_{k_i}-p^{\ell_i}\quad {\text{\rm for}}\quad i=1,2,\ldots,m.
\end{equation}
 We assume $\ell_1>\ell_2>\cdots>\ell_m\ge 0$. Then,
 \begin{equation} \label{4.2}
 F_{k_i}-F_{k_j} = p^{\ell_j}-p^{\ell_j}>0,\quad {\text{\rm for}}\quad 1\le i<j\le m,
 \end{equation}
so $k_i>k_j$. Thus, $k_1>k_2>\cdots >k_m\ge 2$. 

\subsection{$k_1\le 1000$}

Suppose that $k_1\le 1000$. We considered the Diophantine equation $F_{k_2}-F_{k_3}=p^{\ell_2}-p^{\ell_3}$ for $2\le k_3<k_2\le 1000$, $p\ge 5$ and $0\le \ell_3<\ell_2$. Taking $\ell_3=0$, we get 
$F_{k_2}-F_{k_3}+1=p^{\ell_2}$. The above equation has $2161$ solutions $(k_2,k_3,p,\ell_2)$ in the range $1000\ge k_2>k_3\ge 2$ with $\ell_2=1$ and only one solution with $\ell_2>1$ which is $F_{14}-F_{11}+1=17^2$ (other interesting formulas are $F_{12}-F_2+1=12^2$ and $F_{24}-F_{12}+1=215^2$, but $12$ and $215$ are not primes). For each one of these $2162$ quadruples $(k_2,k_3,p,\ell_2)$, we checked whether there exists $k_1\in [k_2+1,1000]$ such that $F_{k_1}-(F_{k_3}-1)=p^{\ell_1}$ for some positive exponent $\ell_1$ and did not find any such instance. This code ran for a few hours in Mathematica. 

Assume next that $\ell_3\ge 1$. We fix $2\le k_3<k_1\le 1000$. Then
$$
(p-1)^{\ell_2}<p^{\ell_2}-p^{\ell_3}=F_{k_2}-F_{k_3}<p^{\ell_2}.
$$
So, $\ell_2\le (\log(F_{k_2}-F_{k_3}))/\log 4$ and once $\ell_2$ is a fixed positive integer in the above range, we have $p=1+\lfloor (F_{k_2}-F_{k_3})^{1/\ell_2}\rfloor$.  Having found $p$, we calculate 
$$
\ell_3=\lfloor\log(p^{\ell_2}-(F_{k_2}-F_{k_3}))/\log p\rfloor,
$$ 
and check whether $\ell_3\ge 1$ and $p^{\ell_3}=p^{\ell_2}-(F_{k_2}-F_{k_3})$. This program ran for a day or so in Mathematica and did not find any solutions. The only solutions found for $F_{k_1}-p^{\ell_1}=F_{k_2}-p^{\ell_2}$ where $1000\ge k_1>k_2\ge 2$ and $\ell_1>\ell_2\ge 1$
were 
$$
F_8-5^2=F_2-5,\quad F_{10}-7^2=F_7-7,\quad F_{12}-11^2=F_{9}-11.
$$
So, our computation shows that there is no integer $c$ having at least three representations as $F_{k_i}-p^{\ell_i}$ with $2\le k_3<k_2<k_1\le 1000$ and some prime $p\ge 5$. 
So, from now on we assume that $k_1>1000$ when $m\ge 3$ and $k_2>1000$ when $m\ge 4$.    

\subsection{Inequalities for $k_i$ in terms of $\ell_i$}

Recall the Binet formula \eqref{eq:Binet} 
$$
F_n = \frac{\alpha^n - \beta^n}{\sqrt{5}}\quad {\text{\rm for all}}\quad n\geq 0.
$$ 
It is well-known and can be easily checked by induction that the inequalities 
 \begin{equation} \label{4.3}
 \alpha^{n -2}\leq F_{n} \leq \alpha^{n -1}\quad {\text{\rm hold for all}}\quad n\geq 1.
 \end{equation}	
Let $i\in \{1,2,\ldots,m-1\}$ and $j\in \{i+1,\ldots,m\}$. Then 
\begin{eqnarray*}
\alpha^{k_i-4} & \le & F_{k_i-2}=F_{k_i}-F_{k_i-1}\le F_{k_i}-F_{k_j}=p^{\ell_i}-p^{\ell_j}<p^{\ell_i}\\
\alpha^{k_i-1} & \ge & F_{k_i}>F_{k_i}-F_{k_j}=p^{\ell_i}-p^{\ell_j}\ge 0.8p^{\ell_i},
\end{eqnarray*}
where for the last inequality we used the fact that $p\ge 5$.  So, we get 
$$
k_i\log \alpha-4\log \alpha<\ell_i\log p<k_i\log \alpha-\log(\alpha/1.25).
$$
Since $4\log \alpha<2$ and $\log(\alpha/1.25)>0.25$, we can record the following lemma.

\begin{lemma}
\label{lem:1}
For $i=1,2,\ldots,m-1$ we have 
$$
k_i\log \alpha-\ell_i\log p\in (\log(\alpha/1.25),4\log \alpha)\subset (c_1,c_2),$$ ${\text{\rm where}}\quad c_1:=0.25,~c_2:=2$.

\end{lemma}

\subsection{Two small linear forms in logarithms}
 
 We assume that $m\ge 3$. 
 We let $(k,\ell):=(k_i,\ell_i)$ for $i=1,2,\ldots,m-2$ and $(k',\ell'):=(k_{j},\ell_{j})$ for some $j=i+1,\ldots,m-1$. Then 
 \begin{equation}
 \label{eq:nn'}
 F_{k}-p^{\ell} = F_{k'}-p^{\ell'}
 \end{equation}
 can be rewritten as
 \begin{eqnarray*}
 \left|\frac{\alpha^k}{\sqrt{5}}-p^\ell\right| & = & \left|\frac{\alpha^{k'}}{\sqrt{5}}-\frac{\beta^{k'}}{\sqrt{5}}+\frac{\beta^k}{\sqrt{5}}-p^{\ell'}\right|\\
 & \le  & \frac{\alpha^{k'}}{\sqrt{5}}+p^{\ell'}+\frac{|\beta|^2+|\beta|^3}{\sqrt{5}}\\
 & = & \frac{\alpha^{k'}}{\sqrt{5}}+p^{\ell'}+\frac{1}{{\sqrt{5}}\alpha}.
 \end{eqnarray*}
 Thus,
 \begin{eqnarray*}
\left|\alpha^k  p^{-\ell} ({\sqrt{5}})^{-1}-1\right| & < & \frac{\alpha^{k'}/{\sqrt{5}}+p^{\ell'}+1/({\sqrt{5}}\alpha)}{p^{\ell}}\\
& < & \frac{(\alpha^4/{\sqrt{5}})p^{\ell'}+p^{\ell'}+1/({\sqrt{5}}\alpha)}{p^{\ell}}\\
& < & \frac{\alpha^4/{\sqrt{5}}+1+1/(p\alpha{\sqrt{5}})}{p^{\ell-\ell'}}<\frac{4.2}{p^{\ell-\ell'}}.
\end{eqnarray*}
 In the above, we used that $\alpha^{k'}<\alpha^4 p^{\ell'}$ which follows from Lemma \ref{lem:1}. 
 So,
 \begin{equation} 
 \label{4.7}
 | \alpha^{k} p^{-\ell} ({\sqrt{5}})^{-1}-1 | <  \frac{4.2}{p^{\ell-\ell'}}.
 \end{equation}	
 We write 
\begin{equation}
\label{eq:Lambdakl}
\Gamma_{k,\ell}:= \alpha^{k} p^{-\ell} ({\sqrt{5}})^{-1}-1.
 \end{equation}
We have that $\Gamma_{k,\ell}\ne 0$, since otherwise $\alpha^{2k}=5p^{2\ell}\in {\mathbb N}$, which is impossible.  
 Inequality \eqref{4.7} shows that 
 \begin{equation}
 \label{eq:Gamma}
 (\ell-\ell')\log p<-\log |\Gamma_{k,\ell}|+\log(4.2),
 \end{equation} 
 and using Lemma \ref{lem:1}, we also have 
\begin{equation}
\label{eq:logalpha}
(k-k')\log \alpha<(\ell-\ell')\log p+(c_2-c_1)<-\log |\Gamma_{k,\ell}|+(1.75+\log(4.2)).
\end{equation}
 This is the first small linear form in logarithms.
We return to equation \eqref{eq:nn'} and use the Binet formula to rewrite it as 
\begin{equation}
\label{eq:2}
\left|\frac{\alpha^k(1-\alpha^{k'-k})}{\sqrt{5}}-p^\ell(1-p^{\ell'-\ell})\right|=\left|\frac{\beta^k}{\sqrt{5}}-\frac{\beta^{k'}}{\sqrt{5}}\right|.
\end{equation}
The above implies that 
\begin{eqnarray*}
\left|\frac{\alpha^k(1-\alpha^{k'-k})}{\sqrt{5}}-p^\ell(1-p^{\ell'-\ell})\right| & \le & \frac{1}{\sqrt{5}}\left(\frac{1}{\alpha^k}+\frac{1}{\alpha^{k'}}\right)\\
& = & \frac{1}{{\sqrt{5}} \alpha^{k'}}\left(1+\frac{1}{\alpha}\right)\\
& = & \frac{\alpha}{{\sqrt{5}} \alpha^{k'}}.
\end{eqnarray*}
Dividing across by $p^\ell(1-p^{\ell'-\ell})$, we get
\begin{equation}
\label{eq:Lambda1}
\begin{aligned}
\left|\alpha^k p^{-\ell} \left(\frac{{\sqrt{5}}(1-p^{\ell'-\ell})}{1-\alpha^{k'-k})}\right)^{-1}-1\right|&<\frac{\alpha}{{\sqrt{5}}\alpha^{k'} p^\ell(1-p^{\ell'-\ell})} \\ &\le \frac{(5/4)(\alpha^4)\alpha}{{\sqrt{5}}\alpha^{k'+k}}<\frac{6.2}{\alpha^{k+k'}},
\end{aligned}
\end{equation}
where we used the  fact that $p\ge 5$ (so $1-p^{\ell'-\ell}\ge 1-1/5)$, as well as the fact that $p^\ell>\alpha^k/\alpha^4$, which follows from Lemma \ref{lem:1}. 
As before, we put
\begin{eqnarray*}
\Gamma_{k,\ell}' & := & \alpha^k p^{-\ell} \left(\frac{{\sqrt{5}}(1-p^{\ell'-\ell})}{1-\alpha^{k'-k})}\right)^{-1}-1;\\
\end{eqnarray*}
Note that $\Gamma_{k,\ell}'\ne 0$ since otherwise \eqref{eq:2} gives that $\beta^k=\beta^{k'}$, so $k=k'$ which is impossible. 
Inequality \eqref{eq:Lambda1} shows that 
$$
(k+k')\log \alpha<-\log |\Gamma_{k,\ell}'|+\log(6.2),
$$
which together with \eqref{eq:logalpha} gives
\begin{equation}
\label{eq:eqkfinal}
k<\frac{1}{2\log \alpha} \left(-\log |\Gamma_{k,\ell}|-\log |\Gamma_{k,\ell}'|+(1.75+\log(4.2)+\log(6.2))\right).
\end{equation}
This is the second small linear form in logarithms. 

\subsection{Bounds on $k$ and $p$}

\begin{lemma}
\label{lem:k1}
If $m\ge 3$, we have:
\begin{itemize}
\item[(i)] $k<7.2\cdot 10^{24} (1+\log k)^2 (\log p)^2$;
\item[(ii)] $k<5\cdot 10^{29} (\log p)^2 (\log\log p)^2$.
\end{itemize}
\end{lemma}

\begin{proof}
We need lower bounds on $\log |\Gamma_{k,\ell}|$ and $\log |\Gamma_{k,\ell}'|$. This we get using Theorem \ref{thm:Mat}. In both cases 
$$
t:=3,~~\alpha_1:=\alpha,~\alpha_2:=p,~~b_1:=k,~b_2:=\ell~~{\text{\rm and}}~~ b_3:=-1.
$$ 
Further, 
$$\alpha_3:={\sqrt{5}}\quad {\text{\rm for}}\quad \Gamma_{k,\ell}\quad {\text{\rm and}}\quad \alpha_3:=\frac{{\sqrt{5}}(1-p^{\ell'-\ell})}{1-\alpha^{k'-k}}\quad {\text{\rm for}}\quad \Gamma_{k,\ell}'.
$$ 
In both cases ${\mathbb K}:={\mathbb Q}(\alpha_1,\alpha_2,\alpha_3)={\mathbb Q}({\sqrt{5}})$ has $D:=2$. 
 Further, we must take 
 $$
 B\ge \max\{|b_1|,|b_2|,|b_3|\}=\max\{k,\ell,1\},
 $$
 and since 
 $$
 \ell \le \frac{k\log \alpha}{\log p}<\frac{k}{3}<k\quad ({\text{\rm because}}\quad p\ge 5>\alpha^3)
 $$
 (see also Lemma \ref{lem:1}), it follows that we can take $B:=k$. Next, we must chose $A_j$ such that 
 $$
 A_j\ge \max\{Dh(\alpha_j),|\log \alpha_j|,0.16\}
 $$
 for $j=1,2,3$. So, we choose 
 $$A_1:=Dh(\alpha_1)=\log \alpha,~A_2:=Dh(\alpha_2)=2\log p
 $$and for $\Gamma_{k,\ell}$ we choose $A_3:=Dh(\alpha_3)=\log 5$. Then, by Theorem \ref{thm:Mat}, 
 we get
 \begin{eqnarray}
 \label{eq:Gammakl}
 \log |\Gamma_{k,\ell}| & > & -1.4\cdot 10^{6}\cdot 3^{4.5}\cdot 2^2(1+\log 2)(1+\log k)\nonumber\\
 & \times & (\log 5)(\log \alpha)(2\log p)\nonumber\\
 & > & -1.51\cdot 10^{12} (\log p)(1+\log k).
 \end{eqnarray}
  Inequalities \eqref{eq:Gamma} and \eqref{eq:logalpha} give 
  \begin{equation}
  \label{eq:max}
  \begin{aligned}
  \max\{(\ell-\ell')\log p,(k-k')\log \alpha\}&<-\log|\Gamma_{k,\ell}|+(1.75+\log(4.2))\\&<1.52 \cdot 10^{12} (\log p)(1+\log k),
  \end{aligned}
  \end{equation}
  so, we can pass to estimate a lower bound for $\Gamma_{k,\ell}'$. We only need to estimate the height of $\alpha_3$:
\begin{eqnarray*}
h(\alpha_3) & \le & h(1-\alpha^{\ell'-\ell})+h(1-p^{k'-k})+h({\sqrt{5}})\\
& \le & h(\alpha^{\ell'-\ell})+h(p^{k'-k})+(1/2)\log 5+2\log 2\\
& \le & (1/2)(\ell-\ell')\log \alpha+(k-k')\log p+(1/2)\log 5+\log 2\\
& < & (1.52/2+1.52)\times 10^{12} (1+\log k)\log p+(1/2)\log 5+2\log 2\\
& < & 2.28 \times 10^{12} (1+\log k)\log p+(1/2)\log 5+2\log 2\\
& < & 2.29\times 10^{12} (1+\log k)\log p,
\end{eqnarray*}
where we used inequality \eqref{eq:max}. 
So, we can take 
$$A_3:=4.6\times 10^{12} (1+\log k)\log p\quad {\text{\rm for}}\quad \Gamma_{k,\ell}'.
$$ 
We get 
\begin{eqnarray*}
\log |\Gamma_{k,\ell}'| & > &  -1.4\cdot 10^6\cdot 3^{4.5} 2^2(1+\log 2)(1+\log k) (\log \alpha)\\
& \times &  (2\log p)(4.6\times 10^{12} (1+\log k) \log p),
\end{eqnarray*}
or simply
\begin{eqnarray}
\label{eq:LowLambda1}
\log |\Gamma_{k,\ell}'| >  -6.91\cdot 10^{24} (1+\log k)^2 (\log p)^2.
 \end{eqnarray}
 Inserting \eqref{eq:LowLambda1} and \eqref{eq:Gammakl} into \eqref{eq:eqkfinal}, we get
 $$
k  <  7.2\cdot 10^{24} (1+\log k)^2(\log p)^2.
 $$
 This is (i).   Assuming $k>10^{10}$, we get 
\begin{eqnarray*}
k & < & 7.2\cdot 10^{24} \left(1+\frac{1}{\log(10^{10})}\right)^2 (\log p)^2 (\log k)^2\\
& < & 7.9\cdot 10^{24}  (\log p)^2(\log k)^2.
\end{eqnarray*}
Finally, we apply Lemma \ref{3.3b} with $s:=2$ and $T:=7.9\cdot 10^{24} (\log p)^2$, to get that 
\begin{eqnarray*}
k & < & 2^2 T (\log T)^2<4\cdot 7.9\cdot 10^{24} (\log p)^2\left(\log(7.9\cdot 10^{24})+2\log\log p\right)^2\\
& < & 31.6 \cdot 10^{24} (\log p)^2 (2\log\log p)^2\left(1+\frac{\log(7.9\cdot 10^{24})}{2\log\log 5}\right)^2\\
& < & 5\cdot 10^{29} (\log p)^2 (\log\log p)^2,
\end{eqnarray*} 
which is (ii). 
\end{proof}

\subsection{An absolute bound on $k_1$}

We assume that $m\ge 4$.
We write inequality \eqref{eq:Lambdakl} in logarithmic form. 
Namely, we put
$$
\Lambda_{k,\ell}:=k\log \alpha-\ell\log p-\log {\sqrt{5}}.
$$
Note that $\Gamma_{k,\ell}=e^{\Lambda_{k,\ell}}-1\ne 0$ so $\Lambda_{k,\ell}\ne 0$. Further, inequality \eqref{eq:Lambdakl} shows that
\begin{equation}
 \label{eq:4.2222}
 |e^{\Lambda_{k,\ell}}-1|<\frac{4.2}{p^{\ell-\ell'}}.
 \end{equation}
 If $\Lambda_{k,\ell}>0$, then 
 $$
 |\Lambda_{k,\ell}|<e^{\Lambda_{k,\ell}}-1<\frac{4.2}{p^{\ell-\ell'}}.
 $$
 If $\Lambda_{k,\ell}<0$, then inequality \eqref{eq:4.2222} together with the fact that $p\ge 5$ implies that 
 $$
 e^{|\Lambda_{k,\ell}|}<\frac{1}{1-\frac{4.2}{5}}=6.25,
 $$
 so 
 $$
 |\Lambda_{k,\ell}|<e^{|\Lambda_{k,\ell}|}\left|1-e^{\Lambda_{k,\ell}}\right|<\frac{6.25\times 4.2}{p^{\ell-\ell'}}=\frac{26.5}{p^{\ell-\ell'}}.
 $$
 Hence, inequality 
 \begin{equation}
 \label{Lambdakl11}
 |\Lambda_{k,\ell}|<\frac{26.5}{p^{\ell-\ell'}}
 \end{equation}
 holds in all cases.  We write inequalities \eqref{Lambdakl11} for 
 $$
 (k,\ell,k',\ell')=(k_i,\ell_i,k_{j},\ell_{j}),~(k_{i+1},\ell_{i+1},k_{j},\ell_{j}),~~{\text{\rm  where}}~~ j\in [i+2,m-1]
 $$ 
 getting 
\begin{eqnarray*}
\left|k_i\log \alpha-\ell_i\log p-\log {\sqrt{5}}\right| & \le & \frac{26.5}{p^{\ell_i-\ell_{j}}},\\
\left|k_{i+1}\log \alpha-\ell_{i+1}\log p-\log {\sqrt{5}}\right| & \le & \frac{26.5}{p^{\ell_{i+1}-\ell_{j}}},
\end{eqnarray*}
and take a linear combination of them to get
\begin{eqnarray}
\label{eq:2form}
|(k_{i+1}\ell_i-k_i\ell_{i+1})\log \alpha-(\ell_i-\ell_{i+1})\log {\sqrt{5}}| & < & \frac{26.5(\ell_i+\ell_{i+1})}{p^{\ell_{i+1}-\ell_{j}}}\nonumber\\
& < & \frac{53\ell_i}{p^{\ell_{i+1}-\ell_{j}}}.
\end{eqnarray}
If the left--hand side is larger than $1/2$, then 
\begin{eqnarray}
\label{eq:31}
p & \le & p^{\ell_{i+1}-\ell_{j}}<106\ell_i\le 106\ell_1<\frac{106 k_1\log \alpha}{\log p}\nonumber\\
& < & \frac{106(\log \alpha)\cdot 5\cdot 10^{29} (\log p)^2(\log\log p)^2}{\log p}\nonumber\\
& < & 3\cdot 10^{31} (\log p)(\log\log p)^2,
\end{eqnarray}
which implies $p<5\cdot 10^{34}$ and next 
$$
k_1<5\cdot 10^{29} (\log p)^2(\log\log p)^2< 7\cdot 10^{34},
$$
which is a pretty good bound on $k_1$. So, assume that the right--hand side of \eqref{eq:2form} is smaller than $1/2$. Then $k_{i+1}\ell_i-k_i\ell_{i+1}$ is positive and 
\begin{eqnarray}
\label{k1k2}
k_{i+1}\ell_i-k_i\ell_{i+1} & < & \frac{(\ell_i-\ell_{i+1})\log {\sqrt{5}}+1/2}{\log\alpha}<\frac{\ell_1\log {\sqrt{5}}}{\log \alpha}\nonumber\\
& < & \frac{k_1\log \alpha\log {\sqrt{5}}}{\log \alpha \log p}\le \frac{k_1}{2}<k_1,
\end{eqnarray}
where we used Lemma \ref{lem:1} and the fact that $p\ge 5$. 
Let $\Lambda$ be the linear form under the absolute value in the left--hand side of \eqref{eq:2form}. It is nonzero since $\alpha$ and ${\sqrt{5}}$ are multiplicatively independent, so if it were zero we would have 
$\ell_i-\ell_{i+1}=0$, which is not the case. Thus, we get
$$
\log p\le (\ell_{i+1}-\ell_{j})\log p<-\log |\Lambda|+\log(53\ell_1).
$$
We need upper bounds on the left--hand side above. The second term has already been estimated in \eqref{eq:31}:
$$
53\ell_1<1.5\cdot 10^{31} (\log p)(\log\log p)^2.
$$ 
As for the first term, we use Theorem \ref{thm:LMN}. We have 
$$
t:=2,~\alpha_1:=\alpha,~\alpha_2:={\sqrt{5}}.
$$ 
We have $D:=2$, $\log A_1:=1/2$, $\log A_2:=(\log 5)/2$. Finally,
$$
b':=\frac{k_{i+1}\ell_i-k_i\ell_{i+1}}{D\log A_2}+\frac{\ell_i-\ell_{i+1}}{D\log A_1}<k_1\left(1+\frac{1}{\log 5}\right)=1.7k_1.
$$
Thus, 
\begin{equation}\nonumber
\begin{aligned}
-\log |\Lambda| &<24.34\cdot D^2 (1/2)((\log5)/2)(\max\{\log b'+0.14, 10.5\})^2\\&<40(\max\{\log b'+0.14,10.5\})^2.
\end{aligned}
\end{equation}
If the maximum is $10.5$, we get 
$$
-\log |\Lambda|<5000.
$$
Thus, in this case 
$$
\log p<5000+\log(1.5\cdot 10^{31} (\log p)^2(\log\log p)^2).
$$
This gives $\log p<5100$. 
If the maximum is in 
\begin{equation}\nonumber
\begin{aligned}
\log b'+0.14=\log(e^{0.14}b')&<\log(e^{0.14}\cdot 1.7 k_1)<\log(2k_1)\\&<\log(10^{30} (\log p)^2(\log\log p)^2),
\end{aligned}
\end{equation}
we get 
$$
-\log |\Lambda|<40 (\log(10^{30}(\log p)^2(\log\log p)^2)^2,
$$
so in this case $\log p$ is smaller than
$$
40 (\log(10^{30}(\log p)^2(\log\log p)^2)^2+\log(1.5\cdot 10^{31} (\log p)^2(\log\log p)^2),
$$
which gives $\log p<4.1\cdot 10^5$. Feeding this into Lemma \ref{lem:k1}, we get 
$$
k_1<5\cdot 10^{29} (\log p)^2(\log\log p)^2<1.5\cdot 10^{43}.
$$
So, we record what we have.

\begin{lemma}
\label{lem:absp}
If $m\ge 4$, we then have $p<e^{4.1\cdot 10^5}$ and $k_1<1.5\cdot 10^{43}$. 
\end{lemma}

\subsection{There are no solutions with $m=4$ and $p<10^{14}$}

The main scope of this section  is to prove the following lemma.

\begin{lemma}
There are no solutions with $m=4$ and $p<10^{14}$. 
\end{lemma}

\begin{proof}
Well, assume that $p<10^{14}$. Lemma \ref{lem:k1} gives
$$
k_1<5\cdot 10^{29} (\log (10^{14})^2 (\log\log 10^{14})^2<10^{34}.
$$
We return to estimate \eqref{eq:2form} with the aim of bounding $p^{\ell_2-\ell_3}$. If the right--had side in \eqref{eq:2form} is at least $1/2$, then 
$$
p^{\ell_2-\ell_3}\le 106 \ell_1<106 k_1<1.1\cdot 10^{36}.
$$
Otherwise, the right--hand side is at most $1/2$, so $k_2\ell_1-k_1\ell_2$ is positive and smaller than $k_1$ as in \eqref{k1k2}. Now $F_{170}>10^{35}>k_1$. We generate the first $171$ convergents of 
$\tau:=\log \alpha/\log {\sqrt{5}}=[0,1,1,\ldots]=[a_0,a_1,a_2,\ldots]$ and get that  $\max\{a_j: 0\le j\le 170\}=330$. Hence, by Lemma \ref{lem:Legendre}, we get that the left--hand side of \eqref{eq:2form} is at least
$$
\frac{1}{(330+2)k_1}.
$$
Thus, we get
$$
p^{\ell_2-\ell_3}<53\cdot 332 \cdot k_1^2<10^{73}.
$$
Next $p^{\ell_3}$ divides $F_{k_i}-F_{k_j}$ for all $i>j\in \{1,2,3\}$. There are two indices $k_i,~k_j$ which are congruent modulo $2$; hence,
$$
F_{k_i}-F_{k_j}=F_{(k_i\pm k_j)/2}L_{(k_i\mp k_j)/2}
$$ 
by Lemma \ref{lem:FL}. 
Let $z(p)$ be the order of appearance of $p$ in the Fibonacci sequence. Since $p< 10^{14}$, Lemma \ref{lem:MR} shows that $p\| F_{z(p)}$. 
Assume that $p^a\| F_{(k_i\pm k_j)/2}$ and $p^b\| L_{(n_i\mp n_j)/2}$. If $a\ge 1$, then $p^{a-1}\mid (k_i\pm k_j)/2$, so $p^{a-1}\le k_1<10^{34}$. Similarly, $p^b\| L_{(k_i\mp k_j)/2}\mid F_{k_i\mp k_j}$, 
so $p^{b-1}\mid (k_i\mp k_j)/2$, so $p^{b-1}<10^{34}$. So, 
$$
p^{\ell_3}\le p^{a-1}\cdot p^{b-1}\cdot p^2<10^{34}\cdot 10^{34} (10^{14})^2<10^{96}.
$$
Thus,
$$
F_{k_2-2}\le F_{k_2}-F_{k_3}<p^{\ell_2}=p^{\ell_2-\ell_3} \cdot p^{\ell_3}<10^{96}\cdot 10^{73}=10^{169},
$$ 
so $k_2<1000$. But we have already shown that in the range $k_2\le 1000$, there are no instances of $k_2>k_3>k_4\ge 2$ and $\ell_2>\ell_3>\ell_4\ge 0$ such that 
$$
F_{k_2}-p^{\ell_2}=F_{k_3}-p^{\ell_3}=F_{k_4}-p^{\ell_4}
$$
with some prime $p\ge 5$. This finishes the proof of the current lemma.
\end{proof}

\subsection{The conclusion}

\begin{lemma}
We have $m\le 4$.
\end{lemma}

\begin{proof}
Assume $m\ge 5$. We return to \eqref{eq:2form} and take $i=1,j=4$. We have 
$$
\ell_1-\ell_4<\ell_1<\frac{k_1\log \alpha}{\log p}<\frac{1.5\cdot 10^{43}\log \alpha}{\log(10^{14})}<2.3\cdot 10^{41}.
$$
Now $2.3\cdot 10^{41}<F_{200}$. We calculated $[a_0,a_1,\ldots,a_{200}]$ for the number $\tau=\log \alpha/\log {\sqrt{5}}$ obtaining $\max\{a_j:0\le j\le 200\}=330$. Hence, the left--hand side of 
\eqref{eq:2form} is at least
$$
\frac{1}{332 k_1},
$$
which gives that 
\begin{equation}
\label{eq:201}
(\ell_2-\ell_4)\log p<\log(332 k_1\ell_1)<\log(332\cdot (2.3\cdot 10^{41}) \cdot (1.5\cdot 10^{43}))<201.
\end{equation}
By Lemma \ref{lem:1}, we get
$$
k_2-k_4<\frac{203}{\log \alpha}<422.
$$
Thus, $k_2-k_3=:a<k_2-k_4=:b$ are in $[1,421]$. Fix $1\le a<b\in [1,421]$. Then $k_2=k_3+a=k_4+b$, so $k_3=k_4+(b-a):=k_4+c$. Hence,
\begin{eqnarray*}
F_{k_4+b}-F_{k_4} & = & p^{\ell_2}-p^{\ell_4}\equiv 0\pmod {p^{\ell_4}};\\
F_{k_4+c}-F_{k_4} & = & p^{\ell_3}-p^{\ell_4}\equiv 0\pmod {p^{\ell_4}}.
\end{eqnarray*}
Using the Binet formula, we get
\begin{eqnarray*}
\alpha^{k_4}(\alpha^b-1) & \equiv & \beta^{k_4}(\beta^b-1)\pmod {p^{\ell_4}}; \\
\alpha^{k_4} (\alpha^c-1) & \equiv & \beta^{k_4}(\beta^c-1)\pmod {p^{\ell_4}}.
\end{eqnarray*}
In the above, for algebraic integers $\gamma,\delta, u$ we write $\gamma\equiv \delta\pmod u$ if $(\gamma-\delta)/u$ is an algebraic integer.  Since $\beta=-\alpha^{-1}$ is a unit, we get
\begin{eqnarray*}
\alpha^{2k_4} (\alpha^b-1) & \equiv & (-1)^{k_4} (\beta^b-1)\pmod {p^{\ell_4}};\\
\alpha^{2k_4} (\alpha^c-1) & \equiv  & (-1)^{k_4} (\beta^c-1)\pmod {p^{\ell_4}}.
\end{eqnarray*}
Multiplying both sides of the second congruence above by $\alpha^b-1$ and using also the first congruence, we get
$$
(-1)^{k_4} (\beta^b-1)(\alpha^c-1)\equiv (-1)^{k_4} (\alpha^b-1)(\beta^c-1)\pmod {p^{\ell_4}}.
$$
Hence, $p^{\ell_4}$ divides
$$
\left|(-1)^c \alpha^{b-c}-\alpha^b-\beta^c+1-((-1)^c \beta^{b-c}-\beta^b-\alpha^c+1)\right|$$
$$=\left|(\alpha^b-\beta^b)-(\alpha^c-\beta^c)\pm (\alpha^{b-c}-\beta^{b-c})\right|.
$$
In particular, 
\begin{equation}
\label{eq:pppp}
p^{\ell_4}\mid F_b-F_c\pm F_{b-c}.
\end{equation}
 We show that $F_b-F_c\pm F_{b-c}$ is nonzero. This is clear if the sign of $F_{b-c}$ is positive since $b>c$. It is also clear if $\max\{c,b-c\}\le b-2$ since then 
$$
F_b-F_{c}-F_{b-c}=F_{b-1}+F_{b-2}-F_c-F_{b-c}>0.
$$
Thus, either $c=b-1$ or $b-c=b-1$. If $c=b-1$, we get 
$$
F_b-F_c-F_{b-c}=F_{b-2}-F_1
$$ 
and this is positive unless $b\in \{2,3,4\}$. Similarly, if $b-c=b-1$, so $c=1$, we get that 
$F_b-F_c-F_{b-c}=F_{b-2}-1$ and again this is positive unless $b\in \{2,3,4\}$. 

If $b=2$, then $a=1$, $k_2=k_3+1=k_4+2$, so $p^{\ell_4}$ divides 
$$
F_{k_2}-F_{k_3}=F_{k_3-1}=F_{k_4}\quad {\text{\rm and~ also}}\quad 
F_{k_2}-F_{k_4}=F_{k_4+1},$$
and this is false since $\gcd(F_{k_4},F_{k_4+1})=1$. 

If $b=3$, then either $a=1$, or $a=2$. When $a=1$, we have $k_2=k_3+1=k_4+3$. So, 
 $$p^{\ell_3}\| F_{k_2}-F_{k_3}=F_{k_3-1}=F_{k_4+1},
 $$
 and also
 $$
p^{\ell_4}\| F_{k_2}-F_{k_4}=F_{k_4+3}-F_{k_4}=F_{k_4+2}+F_{k_4+1}-F_{k_4}=2F_{k_4+1},
$$ 
which implies that $\ell_3=\ell_4$, and this is impossible. If $a=2$, then $k_2=k_3+2=k_4+3$. Thus, 
$$
p^{\ell_4}\| F_{k_4+3}-F_{k_4}=2F_{k_4+1}\quad {\text{\rm and}}\quad p^{\ell_4}\| F_{k_4+1}-F_{k_4}=F_{k_4-1},
$$
and this is impossible since $\gcd(F_{k_4+1},F_{k_4-1})=1$.  

If $b=4$, then $a\in \{1,2,3\}$. If $a=1$, then $k_2=k_3+1=k_4+4$. Then 
$$p^{\ell_3}\| F_{k_2}-F_{k_3}=F_{k_3-1}=F_{k_4+2}
$$
and 
$$p^{\ell_4}\| F_{k_3}-F_{k_4}=F_{k_4+3}-F_{k_4}=2F_{k_4+1},
$$ 
so again 
$p\mid(\gcd(F_{k_4+2},F_{k_4+1})=1$, a contradiction. If $a=2$, then $k_2=k_3+2=k_4+4$, so $p$ divides $F_{k_2}-F_{k_3}=F_{k_3+1}=F_{k_4+3}$ and also $F_{k_3}-F_{k_4}=F_{k_4+1}$. Thus, 
$p$ divides $\gcd(F_{k_4+1},F_{k_4+3})\mid F_{2}=1$, a contradiction. 

Finally, if $a=3$, then $k_2=k_3+3=k_4+4$, so $p$ divides 
$$
F_{k_2}-F_{k_3}=F_{k_2+3}-F_{k_3}=2F_{k_3+1}=2F_{k_4+2}
$$ 
and 
$$F_{k_3}-F_{k_4}=F_{k_4-1}
$$ 
and since $\gcd(F_{k_4+2},F_{k_4-1})\mid F_3=2$, we get a contradiction. 

The above argument shows that the integer which appears in the right--hand side of \eqref{eq:pppp} is nonzero. Its size is at most
$$
F_b+F_{b-c}\le F_{b+1}<\alpha^{421}.
$$
Thus, $p^{\ell_4}<\alpha^{421}$. Since also $p^{\ell_2-\ell_4}<e^{201}$ (see \eqref{eq:201}), we get that 
$$
\alpha^{k_2-4}<F_{k_2-2}\le F_{k_2}-F_{k_3}<p^{\ell_2}=(p^{\ell_2-\ell_4})(p^{\ell_4})<e^{201}\cdot \alpha^{421},
$$
so 
$$
k_2<4+\frac{201+421\log \alpha}{\log \alpha}<850,
$$
but again due to the computation that we did at the beginning, we saw that there do not exist $k_2>k_3>k_4$ in $[1,1000]$ such that $F_{k_2}-p^{\ell_2}=F_{k_3}-p^{\ell_3}=F_{k_4}-p^{\ell_4}$
for some prime $p$ and integers $\ell_2,\ell_3,\ell_4$. This finishes the proof.
\end{proof}

\section*{Acknowledgements} We thank the referees for a careful reading of the first version of this manuscript and for helpful comments and suggestions.

%
%
\end{document}